\documentclass[12pt]{amsart}
\usepackage[margin=1in]{geometry}
\usepackage{amsmath,amssymb,amsthm,mathrsfs,mathtools,enumerate}
\usepackage{enumitem,microtype,xcolor,hyperref}
\hypersetup{colorlinks=true,citecolor=black,filecolor=black,linkcolor=blue,urlcolor=black}
\allowdisplaybreaks
\makeatletter
\@namedef{subjclassname@2020}{\textup{2020} Mathematics Subject Classification}
\makeatother
\newtheorem{theorem}{Theorem}[section]
\newtheorem{lemma}[theorem]{Lemma}
\newtheorem{proposition}[theorem]{Proposition}

\theoremstyle{remark}
\newtheorem{remark}[theorem]{Remark}
\numberwithin{equation}{section}
\newcommand{\R}{\mathbb R}

\newcommand{\one}{\mathbf 1}

\newcommand{\sgn}{\operatorname{sgn}}

\newcommand{\norm}[1]{\lVert#1\rVert}

\title[Maximal operators for Laguerre measures]{Non-centered maximal operators over $\ell^q$ balls for exponential and Laguerre-type measures}
\author[T. A. Bui]{The Anh Bui}
\address{School of Mathematical and Physical Sciences, Macquarie University, NSW 2109, Australia}
\email{the.bui@mq.edu.au}
\author[F. K. Ly]{Fu Ken Ly}
\address{School of Mathematics and Statistics, The Learning Hub, The University of Sydney, Sydney, NSW 2006, Australia}
\email{ken.ly@sydney.edu.au}
\date{}
\begin{document}
\begin{abstract}
In this paper, we study the non-centered Hardy--Littlewood maximal operator over
\(\ell^q\) balls for the measures
\[
 d\mu_\alpha(x)=
 e^{-\sum_{j=1}^d|x_j|}\prod_{j=1}^d|x_j|^{\alpha_j}\,dx,
 \qquad \alpha_j>-1.
\]
We prove that the maximal operator is bounded on \(L^p(d\mu_\alpha)\) for \(1<p<\infty\)
and \(1\le q\le\infty\) in every dimension \(d\ge1\). In particular, for \(\alpha=0\) and \(q=2\), this answers the
question of Nowak, Sasso, Sj\"ogren and Stempak in dimensions \(d\ge5\).
We also prove that the maximal operator is of weak type \((1,1)\) if and only if \(d=1\).
\end{abstract}
\maketitle

\tableofcontents
\section{Introduction}\label{sec:intro}

The Hardy--Littlewood maximal operator is one of the fundamental tools in harmonic analysis and has important applications in PDEs; see, for example, \cite{BrueNguyenStefani}. On a metric measure space \((X,\rho,\mu)\), set
\begin{equation}\label{eq:def-max-general}
 M_{\mu,\rho}f(x)
 =\sup_{\substack{x\in B_\rho\\0<\mu(B_\rho)<\infty}}
 \frac{1}{\mu(B_\rho)}\int_{B_\rho}|f(y)|\,d\mu(y),
\end{equation}
where \(B_\rho\) ranges over open \(\rho\)-balls. Restricting to
\(B_\rho=B_\rho(x,r)\) gives the centered operator. If \(\mu\) is doubling,
both operators are of weak type \((1,1)\) and are bounded on \(L^p(d\mu)\)
for \(1<p<\infty\). On \(\mathbb R^d\), the Besicovitch covering theorem can be used to prove that the centered maximal operator is of weak type \((1,1)\) for Euclidean balls and
locally finite Borel measures. For the non-centered maximal
operator, the weak type \((1,1)\) estimate need not hold in \(d\ge2\).
For example, take the Gaussian measure
\(d\gamma(x)=\pi^{-d/2}e^{-x_1^2-\ldots -x_d^2}\,dx\). With Euclidean balls as
averaging sets, Sjögren \cite{Sjogren83} proved failure of weak type
\((1,1)\) in \(d\ge2\). Forzani, Scotto, Sjögren and Urbina
\cite{Forzani} proved that the same operator is bounded on \(L^p(d\gamma)\)
for all \(1<p<\infty\). Thus the strong \(L^p\) estimates for a
non-doubling measure can hold even if the weak type estimate fails.
It is then natural to ask for which measures and averaging sets this
occurs. For other non-doubling measures we refer to
\cite{AimarForzaniNaibo,LiCuspidalJMPA,SjogrenSoria,DuoNaibo,YangYang2010};
see also \cite{InfanteSoria} for rotationally invariant measures.

For \(x\in\R^d\), write
\[
 |x|_q=\left(\sum_{j=1}^d|x_j|^q\right)^{1/q}\quad(1\le q<\infty),
 \qquad |x|_\infty=\max_{1\le j\le d}|x_j|.
\]
For \(z\in\R^d\) and \(r>0\), let
\(B_q(z,r)=\{x\in\R^d:|x-z|_q<r\}\).
We consider the Laguerre-type measures
\begin{equation}\label{eq:weight}
 d\mu_\alpha(x)=e^{-|x|_1}\prod_{j=1}^d|x_j|^{\alpha_j}\,dx,
 \qquad x\in\R^d,\quad \alpha_j>-1.
\end{equation}
Write \(\mu=\mu_0\) for the exponential measure. Each \(\mu_\alpha\)
is finite and locally doubling: for each fixed \(q\), the doubling
estimate is uniform on balls whose radii are bounded by a fixed
constant, but not on balls of arbitrary radii. The measure \(\mu\) is radial with respect to \(|\cdot|_1\),
and is not Euclidean radial if \(d\ge2\). Write \(M_{\mu,q}\) for
the non-centered maximal operator for \(\mu\) and the balls \(B_q\).

For \(\alpha=0\), Nowak, Sasso, Sjögren and Stempak \cite{NSSS} treated \(q=1,2,\infty\),
that is, diamonds, Euclidean balls, and cubes. For each of these
values of \(q\), \(M_{\mu,q}\) fails to be of weak type \((1,1)\)
if \(d\ge2\). Their \(L^p(d\mu)\) estimates, \(1<p<\infty\), hold
in all dimensions for diamonds and cubes and in \(d\le4\) for
Euclidean balls corresponding to $q=2$. For \(q=2\) they use a rotation to estimate
intersections of balls with triangles when \(d=3\) and with
tetrahedra when \(d=4\). As noted in \cite[Remark~1.1]{NSSS}, the
case \(q=2,\ d\ge5\) was left open.

The restriction of \(\mu_\alpha\) to \(\R^d_+=(0,\infty)^d\) is the
product of the classical Laguerre measures \(x_j^{\alpha_j}e^{-x_j}\,dx_j\).
We write \(M_{\mu_\alpha,q}\) for \eqref{eq:def-max-general} with
measure \(\mu_\alpha\) and balls \(B_q\), and \(q'\) for the conjugate
of \(q\), with \(1/\infty=0\).
The strong \(L^p\) estimates in \cite{NSSS} are for \(\alpha=0\).
Their weak type \((1,1)\) counterexamples for Euclidean balls and
cubes, corresponding to \(q=2\) and \(q=\infty\), respectively, hold for all \(\alpha_j>-1\) when \(d\ge2\). 

In this paper, we prove that \(M_{\mu_\alpha,q}\) is bounded on \(L^p(d\mu_\alpha)\) for every \(1<p<\infty\) and is of weak type \((1,1)\) exactly when \(d=1\), for all \(\alpha\in(-1,\infty)^d\) and \(1\le q\le\infty\). The weak type \((1,1)\) question for diamonds (\(q=1\)) with respect to $\mu_\alpha$ was left open in \cite[Section~3 and Remark~3.1]{NSSS}.

We begin with the strong
\(L^p\) estimate.

\begin{theorem}\label{thm:main}
Let $d\ge1$, $\alpha\in(-1,\infty)^d$, and $1\le q\le\infty$. Then for every $1<p<\infty$ there is a constant $C=C(d,p,q,\alpha)>0$ such that
\begin{equation}\label{eq:main}
 \norm{M_{\mu_\alpha,q}f}_{L^p(\R^d,d\mu_\alpha)}
 \le C\norm{f}_{L^p(\R^d,d\mu_\alpha)}.
\end{equation}
\end{theorem}

We next consider weak type \((1,1)\). More precisely, we have the following result.

\begin{theorem}\label{thm:endpoint}
Let $d\ge1$, $\alpha\in(-1,\infty)^d$ and $1\le q\le\infty$. The
non-centered maximal operator $M_{\mu_\alpha,q}$ is of weak type $(1,1)$
if and only if $d=1$.
\end{theorem}

As in \cite{NSSS}, we study the sections of
\(E=B_q(m,r)\cap\R^d_+\) by the hyperplanes
\(\Sigma(x)=\Sigma(a)+h\), where \(\Sigma(x)=x_1+\cdots+x_d\)
and \(a\) minimizes \(\Sigma\) on \(\overline E\).
The rotation used in \cite{NSSS} for Euclidean balls need not preserve
an $\ell^q$ ball when $q\ne2$. For \(1<q<\infty\), we use the description
of \(a\) in Lemma~\ref{lem:lowest} to estimate the sections in
the original coordinates. For \(q=1\), we choose \(a\) on the
lowest face of the diamond.

Rectangles inside the sections give a lower bound for
\(\mu_\alpha(E)\). To bound the averages over \(E\), we fix some
coordinates, integrate in the remaining coordinates by
\eqref{eq:cartesian-fiber}, and use the $L^p$ estimate for the rectangle
maximal operator in the fixed coordinates.

\medskip

The organization of the paper is as follows. In Section~\ref{sec:reductions}, we reduce the proof of
Theorem~\ref{thm:main} to large truncated balls $E$ in
the positive orthant with $\inf_{x\in E}\Sigma(x)>2$. The proofs of Theorems \ref{thm:main} and \ref{thm:endpoint} will be given in Sections~\ref{sec:geometry} and~\ref{sec:weak-endpoint}{, respectively}.

\section{Preliminaries and reductions}\label{sec:reductions}

\subsection{Notation}

As usual, $q'$ is the conjugate exponent of $q$, where $1'=\infty$
and $\infty'=1$. We use the notation $a\wedge b=\min(a,b)$ and
$a\vee b=\max(a,b)$. The Lebesgue measure of a set $E$ is denoted by
$|E|$.

We write $u_+=\max(u,0)$,  $\one_A$ for the characteristic function of
$A$, and $|I|$ for the cardinality of a finite index set $I$.

For $x\in \R^d$, we write $x_i$ for the $i$-th coordinate of $x$ and set $\Sigma(x)=\sum_{j=1}^d x_j$ and $\R^d_+=(0,\infty)^d$. The {standard} basis for $\R^d$ is denoted by $\{e_1,\ldots, e_d\}$.

The restrictions of $\mu$ and $\mu_\alpha$ to $\R^d_+$, still denoted
by $\mu$ and $\mu_\alpha$, are given by
\begin{equation}\label{eq:mu-plus}
 d\mu(x)=e^{-\Sigma(x)}\,dx,\qquad
 d\mu_\alpha(x)=\prod_{j=1}^d x_j^{\alpha_j}\,d\mu(x),
 \qquad x\in\R^d_+.
\end{equation}

Unless otherwise stated, the constants in $\lesssim$, $\sim$ and
$\simeq$ depend at most on $d,q,\alpha$ and 
on $p$. In particular, they are independent of the center and radius
of the ball. 

\subsection{Reduction by symmetry}

For $x\in\R^d_+$, define
\[
M_{\mu_\alpha,q}^+f(x)
=\sup_{\substack{x\in B_q(m,r)\cap\R^d_+\\
		m\in\overline{\R^d_+},\ r>0}}
\frac{1}{\mu_\alpha(B_q(m,r)\cap\R^d_+)}
\int_{B_q(m,r)\cap\R^d_+}|f(y)|\,d\mu_\alpha(y).
\]

For $x\in\R^d$ and $\varepsilon\in\{-1,1\}^d$, set
\[
 T(x)=(|x_1|,\ldots,|x_d|),\qquad
 R_\varepsilon x=(\varepsilon_1x_1,\ldots,\varepsilon_dx_d).
\]

\begin{lemma}\label{lem:orthant}
Let $1\le q\le\infty$ and $1<p<\infty$. If $M_{\mu_\alpha,q}^+$ is bounded on
$L^p(\R^d_+,d\mu_\alpha)$, then $M_{\mu_\alpha,q}$ is bounded on
$L^p(\R^d,d\mu_\alpha)$.
\end{lemma}

\begin{proof}
For every $c\in\R^d$ and $r>0$,
\begin{equation}\label{eq:fold}
 T(B_q(c,r))\cap\R^d_+
 =B_q(T(c),r)\cap\R^d_+.
\end{equation}
The inclusion 
\[
 T(B_q(c,r))\cap\R^d_+
\subset B_q(T(c),r)\cap\R^d_+
\]
follows from
$\bigl||y_i|-|c_i|\bigr|\le |y_i-c_i|$ for each $i$.
For the reverse inclusion, take $z\in B_q(T(c),r)\cap\R^d_+$
and choose $\varepsilon_i\in\{-1,1\}$ so that
$\varepsilon_i c_i=|c_i|$. Then
$|R_\varepsilon z-c|_q=|z-T(c)|_q<r$ and
$T(R_\varepsilon z)=z$.

Given $f$ on $\R^d$, define
\[
 F(z)=\sum_{\varepsilon\in\{-1,1\}^d}|f(R_\varepsilon z)|,
 \qquad z\in\R^d_+.
\]
The measure is invariant under coordinate reflections. Let $B=B_q(c,r)$ and
$E=B_q(T(c),r)\cap\R^d_+$. Choosing the signs of $c_i$ gives a reflection of $E$ contained in $B$. {By reflecting each orthant part of $B$ into $E$ and using the invariance of $\mu_\alpha$ under these reflections, we obtain}
\[
 \mu_\alpha(E)\le\mu_\alpha(B),\qquad
 \int_B|f|\,d\mu_\alpha\le\int_E F\,d\mu_\alpha .
\]
Hence, for almost every $x\in\R^d$,
\[
 M_{\mu_\alpha,q}f(x)\le M_{\mu_\alpha,q}^+F(T(x)).
\]
Moreover, one has
\[
 \|F\|_{L^p(\R^d_+,d\mu_\alpha)}^p
 \le 2^{d(p-1)}\|f\|_{L^p(\R^d,d\mu_\alpha)}^p, 
\]
which proves the lemma.
\end{proof}

\subsection{Reduction to large balls above a fixed level}

{
We will use the following one-dimensional estimates in the sequel.

\begin{lemma}\label{lem:power}
Let $\beta>-1$. For every $a\ge0$, $L>0$, and fixed $0<\varepsilon\le1$,
\begin{equation}\label{eq:power-one-sided}
 \int_a^{a+\varepsilon L}s^\beta\,ds
 \simeq_{\beta,\varepsilon}L(a+L)^\beta.
\end{equation}
For every $\lambda\ge1$,
\begin{equation}\label{eq:power-enlarge}
 \int_{(a-\lambda L,a+\lambda L)\cap(0,\infty)}s^\beta\,ds
 \le C_\beta\lambda^{1+\beta_+}L(a+L)^\beta,
 \qquad\beta_+=\max(\beta,0).
\end{equation}
Moreover, for $a\ge0$ and $L>0$,
\begin{equation}\label{eq:power-centered}
 \int_{(a-L,a+L)\cap(0,\infty)}s^\beta\,ds
 \simeq_\beta L(a+L)^\beta.
\end{equation}
\end{lemma}
\begin{proof}

By a change of variables,
\[
 \frac{\int_a^{a+\varepsilon L}s^\beta\,ds}{L(a+L)^\beta}
 =\frac{\int_0^\varepsilon(\xi+t)^\beta\,dt}{(1+\xi)^\beta},\qquad \xi=a/L.
\]
The right-hand side is positive and continuous on $[0,\infty)$ and tends to $\varepsilon$ as $\xi\to\infty$. This proves \eqref{eq:power-one-sided}.

{The lower bound in \eqref{eq:power-centered} follows from} \eqref{eq:power-one-sided} with $\varepsilon=1$. If $a\ge2L$, then $s\simeq a+L$ for $s\in(a-L,a+L)$. If $a<2L$, then
\[
 \int_{(a-L)_+}^{a+L}s^\beta\,ds\le\frac{(a+L)^{\beta+1}}{\beta+1}\le\frac{3L(a+L)^\beta}{\beta+1}.
\]
Finally, by \eqref{eq:power-centered} and $\lambda\ge1$,
\[
 \int_{(a-\lambda L,a+\lambda L)\cap(0,\infty)}s^\beta\,ds
 \lesssim_\beta\lambda L(a+\lambda L)^\beta
 \le\lambda^{1+\beta_+}L(a+L)^\beta.
\]

{This completes our proof.}
\end{proof}

\begin{lemma}\label{lem:local-volume}
Let $1\le q\le\infty$, $\alpha_j>-1$, $z\in\overline{\R^d_+}$, and $s>0$. Then
\begin{equation}\label{eq:local-polynomial-mass}
 \int_{B_q(z,s)\cap\R^d_+}\prod_{j=1}^dx_j^{\alpha_j}\,dx
 \simeq_{d,q,\alpha}s^d\prod_{j=1}^d(z_j+s)^{\alpha_j}.
\end{equation}
In particular, for each fixed $R>0$ the measure $\mu_\alpha$ is uniformly doubling on such truncated balls of radii $s\le R$.
\end{lemma}
\begin{proof}

For $c=(2d^{1/q})^{-1}$, we have
\[
 \prod_{j=1}^d(z_j,z_j+cs)\subset B_q(z,s)\cap\R^d_+
 \subset\prod_{j=1}^d((z_j-s)_+,z_j+s).
\]
Hence \eqref{eq:local-polynomial-mass} follows from Lemma~\ref{lem:power}. For $s\le R$ and $x\in B_q(z,2s)\cap\R^d_+$,
\[
 |\Sigma(x)-\Sigma(z)|\le d^{1/q'}|x-z|_q\le2d^{1/q'}R.
\]
Consequently,
\[
 \mu_\alpha(B_q(z,2s)\cap\R^d_+)\simeq_{d,q,\alpha,R}
 e^{-\Sigma(z)}s^d\prod_{j=1}^d(z_j+s)^{\alpha_j}
 \simeq_{d,q,\alpha,R}\mu_\alpha(B_q(z,s)\cap\R^d_+).
\]

{This completes our proof.}
\end{proof}
}

\begin{lemma}\label{lem:small-radii}
	Fix $R>0$.  The part of $M_{\mu_\alpha,q}^+$ obtained by restricting to balls of
	radius at most $R$ is bounded on $L^p(d\mu_\alpha)$ for every $1<p<\infty$.
\end{lemma}
\begin{proof}
	If $x\in E=B_q(m,r)\cap\R^d_+$ and $r\le R$, then
	\[
	E\subset B_q(x,2r)\cap\R^d_+
	\subset B_q(m,3r)\cap\R^d_+.
	\]
	Applying Lemma~\ref{lem:local-volume}  
	(with radius bound $2R$), we have
	\[
	\mu_\alpha(B_q(x,2r)\cap\R^d_+)
	\le \mu_\alpha(B_q(m,4r)\cap\R^d_+)
	\lesssim \mu_\alpha(E).
	\]
	Consequently,
	\[
	\frac{1}{\mu_\alpha(E)}\int_E |f(y)|\,d\mu_\alpha(y)
	\lesssim
	\frac{1}{\mu_\alpha(B_q(x,2r)\cap\R^d_+)}
	\int_{B_q(x,2r)\cap\R^d_+}|f(y)|\,d\mu_\alpha(y)
	\le {M_{\mu_\alpha,q,2R}^c}f(x),
	\]
	{where $M_{\mu_\alpha,q,2R}^c$ is the centered maximal operator
	on $\R^d_+$ restricted to balls $B_q(x,s)\cap\R^d_+$ with $0<s\le2R$. By
	Lemma~\ref{lem:local-volume}, $\mu_\alpha$ is locally doubling,
	so the local centered maximal operator is bounded on $L^p(d\mu_\alpha)$
	for $1<p<\infty$. Taking the supremum over all such balls gives
	the desired result.}
	
	This completes our proof.
	\end{proof}

\bigskip

For a truncated ball $E=B_q(m,r)\cap\R^d_+$, define
\begin{equation}\label{eq:b-def}
 b(E):=\inf_{x\in E}\Sigma(x).
\end{equation}

\begin{lemma}\label{lem:low-level}
Fix $R_0\ge2$.  The maximal operator restricted to balls with $r\ge R_0$
and $b(E)\le2$ is bounded on $L^p(d\mu_\alpha)$ for every $1<p<\infty$.
\end{lemma}

\begin{proof}
Choose $y\in E$ with $\Sigma(y)<9/4$ and put $L=|m-y|_q<r$.  If $L\ge1$, let
\[
 z=y+\frac{m-y}{|m-y|_q};
\]
if $L<1$, put $z=y$.  Since $z\in[y,m]$, we have $z\in\overline{\R^d_+}$. Also $|m-z|_q=L-1$ for $L\ge1$ and $|m-z|_q=L<1$ otherwise. Thus $|m-z|_q+1/2<r$ and
\[
 B_q(z,1/2)\cap\R^d_+\subset E.
\]
Also, we have
$
 \Sigma(z)\le \frac94+d^{1/q'},
$
which implies that, for every $x\in B_q(z,1/2)$,
\[
 \Sigma(x)\le \Sigma(z) +d^{1/q'}\le\frac94+2d^{1/q'}.
\]
{
By Lemma~\ref{lem:local-volume} with $s=1/2$,
\[
 \mu_\alpha(E)\ge \mu_\alpha(B_q(z,1/2)\cap\R^d_+)
 \gtrsim_{d,q,\alpha} \prod_{j=1}^d(z_j+1/2)^{\alpha_j}
 \gtrsim_{d,q,\alpha}1 ,
\]
since $0\le z_j\le 9/4+d^{1/q'}$ for every $j$.
}

Thus,
\[
\frac{1}{\mu_\alpha(E)}\int_E |f(y)|\,d\mu_\alpha(y)\lesssim \|f\|_{L^1(d\mu_\alpha)}.
\] 
 Since $\mu_\alpha(\R^d_+)<\infty$, using H\"older's
inequality,
\[
\|f\|_{L^1(d\mu_\alpha)}\lesssim_p \|f\|_{L^p(d\mu_\alpha)}, \ \ \ 1<p\le \infty.
\]
It follows that
\[
\frac{1}{\mu_\alpha(E)}\int_E |f(y)|\,d\mu_\alpha(y)\lesssim \|f\|_{L^p(d\mu_\alpha)}, \ \ \ \ 1<p<{\infty}.
\] 
{Taking the supremum over these balls gives the desired $L^p$ estimate.}

{This completes our proof.}
\end{proof}

\begin{remark}\label{rem:reductions}
Fix $1<p<\infty$, $1\le q\le\infty$ and $\alpha\in(-1,\infty)^d$.
By Lemma~\ref{lem:orthant}, it suffices to prove the $L^p(d\mu_\alpha)$
estimate of the maximal operator on $\R^d_+$. By Lemmas~\ref{lem:small-radii} and
\ref{lem:low-level}, for any fixed $R_{d,q}\ge2$ we may restrict the
supremum of the maximal operator  to truncated balls satisfying
\begin{equation}\label{eq:admissible}
 E=B_q(m,r)\cap\R^d_+,\qquad m\in\overline{\R^d_+},\qquad
 r\ge R_{d,q},\qquad b(E)>2 .
\end{equation}
We call such balls \emph{admissible} and choose $R_{d,q}$ sufficiently
large in Section~\ref{sec:geometry}.
\end{remark}

\section{\texorpdfstring{Proof of the strong-type $L^p$ estimates}{Proof of the strong-type Lp estimates}}\label{sec:geometry}

When $d=1$, the balls are intervals, irrespective of $q$.
It is {well-known} that the non-centered maximal operator with respect to  locally finite Borel measures
on $\R$ is {of} weak type $(1,1)$, and hence bounded on $L^p({d\mu_\alpha})$ for $1<p<\infty$.

In addition, for $q=\infty$ and $d\ge 1$, let $M_j$ be the non-centered maximal operator in the $j$-th variable with respect to $|x_j|^{\alpha_j}e^{-|x_j|}\,dx_j$. As above, each $M_j$ is bounded on $L^p(d\mu_\alpha)$ for $1<p<\infty$. Since
	\[
	M_{\mu_\alpha,\infty}f(x)\le M_1\cdots M_d f(x),
	\]
	we have, for $1<p<\infty$,
	\[
	\|M_{\mu_\alpha,\infty}f\|_{L^p(d\mu_\alpha)}
	\lesssim_{p,d}\|f\|_{L^p(d\mu_\alpha)}.
	\]
	
Therefore, we may assume $d\ge2$ and $1\le q<\infty$.

\subsection{Strictly convex \texorpdfstring{$\ell^q$}{lq} balls, \texorpdfstring{$1<q<\infty$}{1 < q < infinity}}
Assume $1<q<\infty$. Let $E$ be admissible, i.e., $E$ satisfies \eqref{eq:admissible}, and let
\[
 K:=\overline E=\overline{B_q(m,r)}\cap\overline{\R^d_+}.
\]
Let $a$ minimize $\Sigma$ on the compact set $K$, i.e.,
\[
\Sigma(a) = \min_{x\in K}\Sigma(x).
\]
It is easy to see that $a$ is unique and $a \in \partial B_q(m,r)$. In this case,
\begin{equation}\label{eq:b}
 b(E)=\Sigma(a)=a_1+\cdots+a_d.
\end{equation}
 {When the truncated ball $E$ is fixed, we simply write $b=b(E)$.} 

Define
\begin{equation}\label{eq:IJ}
 \mathcal I(a)=\{i:a_i>0\},
 \qquad \mathcal J(a)=\{1,\ldots,d\}\setminus \mathcal I(a).
\end{equation}
For convenience, we write $\mathcal I$ and $\mathcal J$ when $a$ is fixed. Since $b>2$, the set $\mathcal I$ is nonempty.

\begin{lemma}\label{lem:lowest}
There is a number $\tau>0$ such that
\begin{equation}\label{eq:KKT}
 m_i-a_i=\tau\quad(i\in \mathcal I),
 \qquad 0\le m_j\le\tau\quad(j\in \mathcal J).
\end{equation}
Also,
\begin{equation}\label{eq:tau-equation}
 r^q=|\mathcal I|\tau^q+\sum_{j\in \mathcal J}m_j^q
\end{equation}
and hence
\begin{equation}\label{eq:tau-r}
 d^{-1/q}r\le\tau\le r.
\end{equation}
\end{lemma}

\begin{proof}
As \(0\notin K\), we have \(r<|m|_q\). Set
\[
 g(t)=\sum_{i=1}^d\min\{t,m_i\}^q,\qquad 0\le t\le\max_i m_i.
\]
The function \(g\) is continuous and strictly increasing on this
interval, with \(g(0)=0\) and \(g(\max_i m_i)=|m|_q^q\).
Choose \(\tau\) with \(g(\tau)=r^q\). Then
\[
 a_i^*=(m_i-\tau)_+,\qquad
 |a^*-m|_q^q=\sum_{i=1}^d\min\{\tau,m_i\}^q=r^q,
 \qquad a^*\in K.
\]
For \(x_i\ge0\), the tangent-line inequality gives
\[
 |x_i-m_i|^q
 \ge |a_i^*-m_i|^q
     -q\min\{\tau,m_i\}^{q-1}(x_i-a_i^*)
 \ge |a_i^*-m_i|^q-q\tau^{q-1}(x_i-a_i^*).
\]
In the last inequality, {equality holds} when \(m_i\ge\tau\);
when \(m_i<\tau\), we have \(x_i-a_i^*=x_i\ge0\).
Summing over \(i\), for \(x\in K\) we get
\[
 r^q\ge |x-m|_q^q
 \ge r^q-q\tau^{q-1}\bigl(\Sigma(x)-\Sigma(a^*)\bigr),
 \qquad \Sigma(x)\ge\Sigma(a^*).
\]
Thus \(a^*=a\). The formula for \(a_i^*\) gives \eqref{eq:KKT},
and \(g(\tau)=r^q\) gives \eqref{eq:tau-equation}. Since
\(\mathcal I\ne\emptyset\), we have \(\tau^q\le r^q\le d\tau^q\).

This completes our proof.
\end{proof}

Choose and fix an index
\begin{equation}\label{eq:i0}
 i_0\in \mathcal I,
 \qquad a_{i_0}=\max_{1\le i\le d}a_i.
\end{equation}
If more than one index attains the maximum, we choose the smallest such index.
Then
\begin{equation}\label{eq:large-active}
 a_{i_0}\ge b/d.
\end{equation}

For $c\ge0$ and $u\in\R$, define
\begin{equation}\label{eq:Phi}
 \Phi_c(u)=|c-u|^q-c^q+q c^{q-1}u.
\end{equation}
By the convexity of $s\mapsto |s|^q$,
\begin{equation}\label{eq:Phi-positive}
 \Phi_c(u)\ge0 \ \ \text{for all $u\in \R$}.
\end{equation}
For $j\in \mathcal J$, {by} \eqref{eq:KKT}, {we have}
\begin{equation}\label{eq:sj}
 s_j=\tau^{q-1}-m_j^{q-1}\ge0.
\end{equation}

\begin{lemma}\label{lem:cap}
Let $x=a+y\in E$ and assume
\begin{equation}\label{eq:height-h}
h:=\Sigma(y)\ge0.
\end{equation}
Then $y_j\ge0$ for $j\in \mathcal J$ and
\begin{equation}\label{eq:cap}
 \sum_{i\in \mathcal I}\Phi_\tau(y_i)
 +\sum_{j\in \mathcal J}\Phi_{m_j}(y_j)
 +q\sum_{j\in \mathcal J}s_jy_j
 \le q\tau^{q-1}h.
\end{equation}
\end{lemma}

\begin{proof}
If $j\in \mathcal J$, then $a_j=0$. Hence, $y_j=x_j\ge0$ for all $j\in \mathcal J$.  Since
$x\in B_q(m,r)$ and $a\in\partial B_q(m,r)$, we have
$|x-m|_q\le |a-m|_q$. It follows that
\[
 \sum_i|a_i+y_i-m_i|^q
 \le \sum_i|a_i-m_i|^q.
\]
Since $a_i-m_i=-\tau$ for $i\in \mathcal I$ and $a_j=0$ for $j\in \mathcal J$, the preceding inequality implies
\[
 \sum_{i\in \mathcal I}|y_i-\tau|^q+ \sum_{j\in \mathcal J}|y_j-m_j|^q
\le \sum_{i\in \mathcal I}\tau^q +\sum_{j\in \mathcal J}|m_j|^q.
\]
Hence,
\[
\sum_{i\in \mathcal I}\Big[|y_i-\tau|^q-\tau^q\Big]+ \sum_{j\in \mathcal J}\Big[|y_j-m_j|^q-|m_j|^q\Big]
\le 0.
\]
This, together with \eqref{eq:Phi}, implies 
\begin{align*}
 0&\ge
 \sum_{i\in \mathcal I}\big(\Phi_\tau(y_i)-q\tau^{q-1}y_i\big)
 +\sum_{j\in \mathcal J}\big(\Phi_{m_j}(y_j)-qm_j^{q-1}y_j\big)\\
 &=\sum_{i\in \mathcal I}\Phi_\tau(y_i)
 +\sum_{j\in \mathcal J}\Phi_{m_j}(y_j)
 -q\tau^{q-1}h+q\sum_{j\in \mathcal J}s_jy_j.
\end{align*}
This yields \eqref{eq:cap}.

{This completes our proof.}
\end{proof}

We need two elementary consequences of one-dimensional convexity.

\begin{lemma}\label{lem:scalar}
Fix $1<q<\infty$.  There are constants $\eta_q\in(0,1)$ and $C_q>0$
with the following properties.

\begin{enumerate}[label=\textup{(\roman*)}]
\item Let $\tau>0$, $0<h\le\eta_q\tau$ and suppose that
$
 \Phi_\tau(u)\le q\tau^{q-1}h.
$
Then
\begin{equation}\label{eq:scalar-active}
 |u|\le C_q\sqrt{\tau h}.
\end{equation}

\item Let {$\tau>0$,} $0\le c\le\tau$, $s=\tau^{q-1}-c^{q-1}$, and $u\ge0$.  If
$0<h\le\eta_q\tau$ and
$
 \Phi_c(u)+qsu\le q\tau^{q-1}h,
$
then
\begin{equation}\label{eq:scalar-inactive}
 u\le C_q\sqrt{\tau h},
 \qquad
 su\le\tau^{q-1}h.
\end{equation}
\end{enumerate}
\end{lemma}

\begin{proof}
For (i), write $u=\tau v$ and $h=\tau\delta$.  {Then the condition $\Phi_\tau(u)\le q\tau^{q-1}h$ becomes}
\[
f(v):= |1-v|^q-1+qv\le q\delta.
\]
By convexity and \(f(0)=0\), choose \(\eta_q\) with
\[
 q\eta_q<\min\{f(-1/2),f(1/2)\}.
\]
Then \(f(v)\le q\delta\le q\eta_q\) implies \(|v|<1/2\). Since
\(f(0)=f'(0)=0\) and \(f''(v)=q(q-1)(1-v)^{q-2}\), we have
\[
 v^2\lesssim_q f(v)\le q\delta,\qquad
 u^2=\tau^2v^2\lesssim_q \tau h.
\]

For (ii), \(0\le qsu\le\Phi_c(u)+qsu\le q\tau^{q-1}h\).
Fix \(\kappa_q=1/2\). If \(s\ge\kappa_q\tau^{q-1}\), then
\(u\le 2h\le 2\sqrt{\tau h}\).
If $s<\kappa_q\tau^{q-1}$, then $c\ge c_q\tau$.  After decreasing
$\eta_q$, we may apply the argument of (i), now at scale $c$, to
$\Phi_c(u)\lesssim_q c^{q-1}h$.  It gives
$u\lesssim_q\sqrt{ch}\le\sqrt{\tau h}$.

{This completes our proof.}
\end{proof}

\begin{lemma}\label{lem:upper-Taylor}
Fix $1<q<\infty$.  There is $C_q$ such that for $c\ge0$ and $u\in\R$,
\begin{equation}\label{eq:Taylor-general}
 \Phi_c(u)\le
 \begin{cases}
 C_q\big(c^{q-2}u^2+|u|^q\big),&q\ge2,\\
 C_q|u|^q,&1<q<2.
 \end{cases}
\end{equation}
Moreover, let $c_0\in(0,1]$.  If
\begin{equation}\label{eq:Taylor-comparable-assumptions}
 c_0\tau\le c\le\tau,\qquad \tau\ge1,
 \qquad |u|\le\varepsilon\sqrt\tau,
\end{equation}
where $0<\varepsilon\le c_0/2$, then
\begin{equation}\label{eq:Taylor-comparable}
 \Phi_c(u)\le C_{q,c_0}\varepsilon^2\tau^{q-1}.
\end{equation}
\end{lemma}

\begin{proof}
Let
\[
 F(t)=|t|^q,\qquad t\in\R.
\]
Since $c\ge0$,
\begin{equation}\label{eq:Phi-as-remainder}
 \Phi_c(u)=F(c-u)-F(c)+F'(c)u.
\end{equation}
We first prove \eqref{eq:Taylor-general}, with \(c^0=1\) when \(q=2\).

Assume $q\ge2$.  Then $F\in C^2(\R)$ and
\[
 |F''(t)|\le C_q|t|^{q-2}.
\]
Taylor's formula at $c$ gives
\begin{equation}
	\label{eq:Taylor-integral}\begin{aligned}
 \Phi_c(u)
 &=u^2\int_0^1(1-s)F''(c-su)\,ds\\
 &\le C_qu^2\int_0^1|c-su|^{q-2}\,ds\\
 &\le C_q(c+|u|)^{q-2}u^2
 \le C_q\big(c^{q-2}u^2+|u|^q\big).
\end{aligned}
\end{equation}
This proves the first line of \eqref{eq:Taylor-general}.

Now assume $1<q<2$.  Since 
\[
 F'(t)=q\,\sgn(t)|t|^{q-1},\qquad \sgn(0)=0,
\]
we have
\begin{equation}\label{eq:Fprime-holder}
 |F'(s)-F'(t)|\le C_q|s-t|^{q-1},
 \qquad s,t\in\R.
\end{equation}
Using \eqref{eq:Phi-as-remainder} and the fundamental theorem of calculus,
\begin{align*}
 \Phi_c(u)
 &=u\int_0^1\big[F'(c)-F'(c-su)\big] \,ds.
\end{align*}
Hence \eqref{eq:Fprime-holder} yields
\[
 0\le\Phi_c(u)
 \le C_q|u|\int_0^1(s|u|)^{q-1}\,ds
 \le C_q|u|^q.
\]
This proves the second inequality of \eqref{eq:Taylor-general}.

It remains to prove \eqref{eq:Taylor-comparable}.  Suppose
\eqref{eq:Taylor-comparable-assumptions} holds.  Since $\tau\ge1$ and
$\varepsilon\le c_0/2$,
\[
 |u|\le\varepsilon\sqrt\tau\le\varepsilon\tau
 \le\frac{c_0\tau}{2}\le\frac c2.
\]
Consequently, for every $0\le s\le1$,
\begin{equation}\label{eq:segment-comparable}
 \frac c2\le c-su\le\frac{3c}{2}.
\end{equation}
{Since the segment in \eqref{eq:segment-comparable} stays away from the origin, the second-order integral remainder formula is valid for every $1<q<\infty$.  Hence}
\begin{equation}\label{eq:Phi-integral-second}
 {
 \Phi_c(u)
 =u^2\int_0^1(1-s)F''(c-su)\,ds
 \le C_q u^2 c^{q-2}.
 }
\end{equation}
{Since $c_0\tau\le c\le\tau$ and $|u|\le \varepsilon\sqrt\tau$, it follows that, for every $1<q<\infty$,}
\[
{
\Phi_c(u)\le C_{q,c_0}\varepsilon^2 \tau^{q-1}.
}
\]

{This completes our proof.}
\end{proof}

Let $r\ge R_{d,q}$ be large enough so that $\tau\ge1$, where $R_{d,q}$ will be fixed later. Define, for each $j\ne i_0$, the length
\begin{equation}\label{eq:Lj-strict}
 L_j=\begin{cases}
 \min\{b+1,\sqrt r\},&j\in \mathcal I\setminus\{i_0\},\\
 \min\{b+1,\sqrt r,\tau^{q-1}/s_j\},&j\in \mathcal J,\ s_j>0,\\
 \min\{b+1,\sqrt r\},&j\in \mathcal J,\ s_j=0.
 \end{cases}
\end{equation}
The bounds \eqref{eq:tau-r} and $s_j\le\tau^{q-1}$ show that
\begin{equation}\label{eq:Lj-range}
 1\le L_j\le b+1.
\end{equation}

\begin{lemma}\label{lem:strict-geometry}
Let $1<q<\infty$. There exist $R_{d,q}>0$, $\varepsilon_{d,q}>0$, and $C_{d,q}<\infty$ such that, for every admissible $E$ with $r\ge R_{d,q}$, the following hold, with \(c_d=1/(2d)\).
\begin{enumerate}[label=\textup{(\roman*)}]
\item For every $x\in E$ with $h=\Sigma(x)-b>0$,
\begin{equation}\label{eq:strict-outer}
 |x_j-a_j|\le C_{d,q}(1+h)L_j,
 \qquad j\ne i_0.
\end{equation}
\item For every $1<h<2$, the following box in the hyperplane $\Sigma(x)=b+h$ given by
\begin{equation}\label{eq:strict-inner}
 Q_h=\left\{a+h e_{i_0}+\sum_{j\ne i_0}z_j(e_j-e_{i_0}):
                 0<z_j<\varepsilon_{d,q}L_j\right\},
\end{equation}
is contained in $E$. On this {box},
\begin{equation}\label{eq:distinguished-comparable}
 c_db\le x_{i_0}\le 2b.
\end{equation}
\end{enumerate}
\end{lemma}
\begin{proof}

(i) Set $y=x-a$. By Lemmas~\ref{lem:cap} and \ref{lem:scalar}, if $0<h\le\eta_q\tau$, then
\[
 |y_j|\lesssim_{d,q}\sqrt{rh}\quad(j\in \mathcal I),\qquad
 {0\le y_j\lesssim_{d,q}\sqrt{rh}\quad(j\in \mathcal J).}
\]
Also $|y_j|\le b+h$ since $0\le x_j,a_j\le b+h$.
For $j\in\mathcal J$ with $s_j>0$, the  inequality
\eqref{eq:cap} gives $s_jy_j\le\tau^{q-1}h$. Therefore,
\[
 0\le y_j\lesssim_{d,q}(1+h)
 \min\{b+1,\sqrt r,\tau^{q-1}/s_j\}=(1+h)L_j.
\]
For $j\in\mathcal I\setminus\{i_0\}$
or $j\in\mathcal J$ with $s_j=0$, the first two bounds give
$$|y_j|\lesssim_{d,q}(1+h)\min\{b+1,\sqrt r\}=(1+h)L_j.
$$
If $h>\eta_q\tau$, then
\[
 |y_j|\le |x_j-m_j|+|m_j-a_j|\le2r\lesssim_{d,q}h\le(1+h)L_j.
\]

(ii) Write $\varepsilon=\varepsilon_{d,q}$ and $S=\sum_{j\ne i_0}z_j$. By \eqref{eq:Lj-strict} and \eqref{eq:large-active},
\[
 S\le(d-1)\varepsilon\min\{b+1,\sqrt r\},\qquad S\le a_{i_0}/2
\]
for $\varepsilon$ sufficiently small. Hence,
\[
 b/(2d)\le a_{i_0}/2\le x_{i_0}\le b+h\le2b.
\]
For \(x=a+y\in Q_h\),
\[
 y_j=z_j\ (j\ne i_0),\qquad y_{i_0}=h-S,\qquad \Sigma(y)=h.
\]
By \eqref{eq:Lj-strict} and \eqref{eq:Taylor-comparable},
\begin{equation}\label{eq:strict-proof-active}
 \Phi_\tau(z_j)\lesssim_{d,q}\varepsilon^2\tau^{q-1},\qquad j\in \mathcal I\setminus\{i_0\}.
\end{equation}
Since $|h-S|\le{2+C_{d,q}\varepsilon\sqrt\tau}\le\tau/2$ for $\tau$ large, Taylor's formula gives
\begin{equation}\label{eq:strict-proof-i0}
 \Phi_\tau(h-S)\lesssim_q\tau^{q-2}(4+S^2)
 \lesssim_{d,q}(\tau^{-1}+\varepsilon^2)\tau^{q-1}.
\end{equation}
For $j\in\mathcal J$ with $s_j>0$, the second line of
\eqref{eq:Lj-strict} and the definition of $Q_h$ imply
\[
 0<s_jz_j<\varepsilon s_j L_j\le\varepsilon\tau^{q-1}.
\]
When $s_j=0$, the term $s_jz_j$ vanishes. Put $\theta_q=1/2$. If $s_j\ge\theta_q\tau^{q-1}$, then $z_j\le\varepsilon/\theta_q$; by \eqref{eq:Taylor-general},
\[
 \Phi_{m_j}(z_j)\lesssim_q\varepsilon^{\min\{2,q\}}\tau^{q-1}.
\]
If $s_j<\theta_q\tau^{q-1}$, then $m_j\simeq_q\tau$ and $z_j\lesssim_{d,q}\varepsilon\sqrt\tau$. \eqref{eq:Taylor-comparable} gives $\Phi_{m_j}(z_j)\lesssim_{d,q}\varepsilon^2\tau^{q-1}$. Therefore
\begin{align*}
 &\sum_{i\in \mathcal I}\Phi_\tau(y_i)+\sum_{j\in \mathcal J}\Phi_{m_j}(y_j)+q\sum_{j\in \mathcal J}s_jy_j\\
 &\qquad\le C_{d,q}(\varepsilon+\varepsilon^{\min\{2,q\}}+\tau^{-1})\tau^{q-1}
 <q\tau^{q-1}h.
\end{align*}
Here we choose \(\varepsilon\) small and then \(R_{d,q}\) large. Subtracting
\(q\tau^{q-1}h\) as in Lemma~\ref{lem:cap} gives
\[
 |x-m|_q^q-r^q<0,\qquad Q_h\subset E .
\]
This completes our proof.

\end{proof}

\subsection{Diamonds, \texorpdfstring{$q=1$}{q = 1}}

We now prove the corresponding estimates for $q=1$. Here \(R_d=R_{d,1}\) is
a sufficiently large constant depending only on \(d\), and
\(c_d=1/(2d)\).

\begin{lemma}\label{lem:diamond-geometry}
Let $E=B_1(m,r)\cap\R^d_+$ {with $r\ge R_d\ge1$ and $b(E)>2$}. Put
\begin{equation}\label{eq:diamond-choice}
 M=\Sigma(m),\quad b=M-r,
 \quad a=\frac bM m,\quad \delta_j=m_j-a_j=\frac rM m_j.
\end{equation}
{Choose $i_0=\min\{i:m_i=\max_j m_j\}$.} Then $a\in\overline E$ minimizes $\Sigma$, and
\begin{equation}\label{eq:diamond-largest}
 a_{i_0}\ge b/d,\qquad\delta_{i_0}\ge r/d.
\end{equation}
For $j\ne i_0$, set
\begin{equation}\label{eq:diamond-Lj}
 L_j=\min\{m_j+1,r,b+1\}.
\end{equation}
Then $1\le L_j\le b+1$, and there exist $\varepsilon_d>0$ and $C_d<\infty$ such that:
\begin{enumerate}[label=\textup{(\roman*)}]
\item Every $x\in E$ with $h=\Sigma(x)-b>0$ satisfies
\begin{equation}\label{eq:diamond-outer}
 |x_j-a_j|\le C_d(1+h)L_j\quad(j\ne i_0).
\end{equation}
\item For each $j\ne i_0$ choose
\begin{equation}\label{eq:diamond-inner-interval}
 I_j=\begin{cases}
 (a_j-2\varepsilon_dL_j,a_j-\varepsilon_dL_j),&a_j\ge L_j/2,\\
 (a_j+\varepsilon_dL_j,a_j+2\varepsilon_dL_j),&a_j<L_j/2.
 \end{cases}
\end{equation}
For every $1<h<2$, the box defined by $x_j\in I_j$ for $j\ne i_0$ and $x_{i_0}=b+h-\sum_{j\ne i_0}x_j$ is contained in $E$. On this box, $c_db\le x_{i_0}\le2b$; furthermore, for every $\beta>-1$,
\begin{equation}\label{eq:diamond-interval-mass}
 \int_{I_j}s^\beta\,ds\gtrsim_{\beta,d}
 L_j(a_j+L_j)^\beta.
\end{equation}
\end{enumerate}
\end{lemma}
\begin{proof}

Since $\Sigma(a)=M-r$ and $|m-a|_1=r$, we have $a\in\overline E$ and
\[
 \Sigma(x)\ge M-|m-x|_1>M-r=\Sigma(a),\qquad x\in E.
\]
Also, $m_{i_0}\ge M/d$ implies \eqref{eq:diamond-largest}.

(ii) Choose \(\varepsilon=\varepsilon_d\in(0,1/8)\) sufficiently small and
take the intervals \(I_j\) in \eqref{eq:diamond-inner-interval}. In both cases,
\[
 \begin{gathered}
 |I_j|=\varepsilon L_j,\qquad s\simeq_d a_j+L_j\quad(s\in I_j),\\
 \int_{I_j}s^\beta\,ds\simeq_{\beta,d} L_j(a_j+L_j)^\beta\quad(\beta>-1).
 \end{gathered}
\]
Fix $1<h<2$ and $x_j\in I_j$ for $j\ne i_0$.  Reconstruct
	$x_{i_0}$ by $\Sigma(x)=b+h$, and write
	$u_j=x_j-a_j$, $u_{i_0}=h-\sum_{j\ne i_0}u_j$.
	From \eqref{eq:diamond-largest}, \eqref{eq:diamond-Lj} and
	$L_j\le b+1$, choosing $\varepsilon$ sufficiently small gives
	\[
		x_{i_0}\ge a_{i_0}/2\ge b/(2d).
	\]
	Furthermore, $|\sum_{j\ne i_0}u_j|\le 2\varepsilon(d-1)r$.
	Since $\delta_{i_0}\ge r/d$, increasing $R_d$ if needed gives
	$x_{i_0}<m_{i_0}$.
	
	We claim that $\sum_j(x_j-m_j)_+<h/2$.
	Its $i_0$ term is zero.  If $a_j\ge L_j/2$, then
	$x_j<a_j\le m_j$.  Otherwise, $u_j\le2\varepsilon L_j$.
	For $L_j\ge4$, the inequality $L_j\le m_j+1$
	implies
	\[
	\delta_j=m_j-a_j\ge L_j-1-a_j
	>L_j/2-1\ge L_j/4.
	\]
	Thus $x_j\le m_j$ when $\varepsilon<1/8$.  For $L_j<4$,
	one has $(x_j-m_j)_+\le8\varepsilon$.
	Taking $8d\varepsilon<1/2$, we obtain the claim.
	The identity
	\begin{equation}\label{eq:diamond-excess}
		|x-m|_1=r-h+2\sum_j(x_j-m_j)_+
	\end{equation}
	now shows that $x\in E$.
All coordinates of \(x\) are positive by construction, so
\(x_{i_0}\le\Sigma(x)=b+h<2b\). This proves part (ii).

(i) Let $x\in E$ satisfy $\Sigma(x)=b+h$, where $h>0$.
	By \eqref{eq:diamond-excess},
	\[
	\sum_{k=1}^d(x_k-m_k)_+<\frac h2,
	\]
	and hence $x_j<m_j+h/2$ for every $j$.  Fix $j\ne i_0$.
	If $x_j\le a_j$, then $x_j\ge0$, $\Sigma(a)=b$, and
	$|x_j-m_j|<r$ imply
	\[
	 a_j-x_j\le \min\{m_j,b,r\}\le L_j.
	\]
	If $x_j>a_j$, we use $x_j<m_j+h/2$, $x_j\le b+h$,
	and $|x_j-m_j|<r$ to obtain
	\[
	\begin{aligned}
	 x_j-a_j&\le m_j+\tfrac h2\le(1+h)(m_j+1),\\
	 x_j-a_j&\le b+h\le(1+h)(b+1),\\
	 x_j-a_j&<r+(m_j-a_j)\le2r.
	\end{aligned}
	\]
	Consequently, by the definition of $L_j$,
	\[
	 |x_j-a_j|\le2(1+h)L_j,\qquad j\ne i_0.
	\]
This proves part (i).
\end{proof}

\subsection{Weighted section estimates}\label{sec:weighted-geometry}
Fix $1\le q<\infty$ and retain the notation $E,b,a,i_0$ and $L_j$ from the preceding two subsections.

\begin{proposition}\label{prop:mass-lower}
{For every $E$ as in \eqref{eq:admissible},}
\begin{equation}\label{eq:mass-lower}
 \mu_\alpha(E)\gtrsim_{d,q,\alpha}
 e^{-b}b^{\alpha_{i_0}}\prod_{j\ne i_0}L_j(a_j+L_j)^{\alpha_j}.
\end{equation}
\end{proposition}
\begin{proof}
{For $j\ne i_0$, let $I_j=(a_j,a_j+\varepsilon_{d,q}L_j)$ if $q>1$, and let $I_j$ be as in \eqref{eq:diamond-inner-interval} if $q=1$. By Lemmas~\ref{lem:strict-geometry} and \ref{lem:diamond-geometry}, for $1<h<2$ and $x_j\in I_j$ ($j\ne i_0$), we have}
\[
 \Sigma(x)=b+h,\qquad x_{i_0}=b+h-\sum_{j\ne i_0}x_j\simeq_d b,
\]
{and $x\in E$.} {For $t=\Sigma(x)$ and $\widehat x=(x_j)_{j\ne i_0}$, we have $dx=dt\,d\widehat x$.} Hence
\begin{align*}
 \mu_\alpha(E)
 &\ge\int_1^2e^{-(b+h)}
 \int_{\prod_{j\ne i_0}I_j}x_{i_0}^{\alpha_{i_0}}
 \prod_{j\ne i_0}x_j^{\alpha_j}\,d\widehat x\,dh\\
 &\gtrsim_{d,\alpha}e^{-b}b^{\alpha_{i_0}}
 \prod_{j\ne i_0}\int_{I_j}x_j^{\alpha_j}\,dx_j\\
 &\gtrsim_{d,q,\alpha}e^{-b}b^{\alpha_{i_0}}
 \prod_{j\ne i_0}L_j(a_j+L_j)^{\alpha_j}.
\end{align*}

Here we used Lemma~\ref{lem:power} for $q>1$ and \eqref{eq:diamond-interval-mass} for $q=1$.

{This completes our proof.}
\end{proof}

\subsection{Marginals in the original coordinates}\label{sec:cartesian-marginal}

\leavevmode{Fix an index $i_0\in\{1,\ldots,d\}$ and let $J$ and $K$ be any index sets}
\[
 J\subseteq\{1,\ldots,d\}\setminus\{i_0\},\qquad
 K=\{1,\ldots,d\}\setminus(J\cup\{i_0\}).
\]
We retain the \(J\)-coordinates and integrate out the \(K\)-coordinates.

For \(x\in\R^d_+\), write $v=(x_j)_{j\in J}$, $y=(x_k)_{k\in K}$ and
\[
 t=\Sigma(x),\qquad s=t-\sum_{j\in J}v_j,\qquad
 d\omega_J(v)=\prod_{j\in J} v_j^{\alpha_j}\,dv_j ,\qquad
 \beta=\alpha_{i_0}+\sum_{k\in K}(\alpha_k+1).
\]

For each fixed $t,v$ with $s>0$, the remaining coordinates satisfy
$\sum_{k\in K}y_k<s$ and $x_{i_0}=s-\sum_{k\in K}y_k$. Since the
Jacobian of $x\mapsto(t,v,y)$ has absolute value one, integration in
these coordinates reduces to
\begin{equation}\label{eq:cartesian-fiber}
 \int_{\substack{y_k>0\\ \sum_{k\in K}y_k<s}}
 \left(s-\sum_{k\in K}y_k\right)^{\alpha_{i_0}}
 \prod_{k\in K}y_k^{\alpha_k}\,dy=C_{\alpha,J,i_0}s^\beta .
\end{equation}
To prove \eqref{eq:cartesian-fiber} when $K\ne\emptyset$, fix $k\in K$ and put
$r=s-\sum_{\ell\in K\setminus\{k\}}y_\ell$. The integral in $y_k$ is
\[
 \int_0^r(r-y_k)^{\alpha_{i_0}}y_k^{\alpha_k}\,dy_k
 =B(\alpha_k+1,\alpha_{i_0}+1)\,r^{\alpha_{i_0}+\alpha_k+1}.
\]
Here \(B\) is the beta function. Repeating the calculation and writing
\(\Gamma\) for the gamma function gives
\[
 C_{\alpha,J,i_0}
 =\frac{\Gamma(\alpha_{i_0}+1)\prod_{k\in K}\Gamma(\alpha_k+1)}
        {\Gamma(\beta+1)}.
\]

Set $c_d=1/(2d)$. For $s\ge c_dt$ we have
\[
 c_d\le \frac{s}{t}\le 1,\qquad
 e^{-t}C_{\alpha,J,i_0}s^\beta
 = C_{\alpha,J,i_0}\left(\frac{s}{t}\right)^\beta e^{-t}t^\beta
 \simeq_{\alpha,d} e^{-t}t^\beta.
\]
Fix \(1<p<\infty\) and \(f\in L^p(d\mu_\alpha)\). For \(t>0\)
and \(v\in(0,\infty)^J\), write $x(t,v,y)$ for the point with coordinates
$x_j=v_j$ ($j\in J$), $x_k=y_k$ ($k\in K$) and
$x_{i_0}=s-\sum_{k\in K}y_k$ whenever \(s>0\) and
\(y_k>0,\ \sum_{k\in K}y_k<s\). Define
\[
 F_J(t,v)=
 \begin{cases}
 \displaystyle
 \frac{\int_{\substack{y_k>0\\\sum_{k\in K}y_k<s}}
 |f(x(t,v,y))|
 \left(s-\sum_{k\in K}y_k\right)^{\alpha_{i_0}}
 \prod_{k\in K}y_k^{\alpha_k}\,dy}
 {C_{\alpha,J,i_0}s^\beta},& s\ge c_dt,\\[1.3em]
 0,& s<c_dt.
 \end{cases}
\]
In particular, if $K=\emptyset$ and $s\ge c_dt$, let $x(t,v)$ have
$x_j=v_j$ for $j\in J$ and $x_{i_0}=s$. Then
\[
 F_J(t,v)=\frac{|f(x(t,v))|s^{\alpha_{i_0}}}{s^{\alpha_{i_0}}}
          =|f(x(t,v))|.
\]
For $s\ge c_dt$, set
\[
 D_s=\left\{y\in(0,\infty)^K:\sum_{k\in K}y_k<s\right\},
 \qquad
 w_s(y)=\left(s-\sum_{k\in K}y_k\right)^{\alpha_{i_0}}
             \prod_{k\in K}y_k^{\alpha_k}.
\]
By H\"older's inequality and \eqref{eq:cartesian-fiber},
\begin{align*}
 \left(\int_{D_s}|f(x(t,v,y))|w_s(y)\,dy\right)^p
 &\le \left(\int_{D_s}|f(x(t,v,y))|^p w_s(y)\,dy\right)
       \left(\int_{D_s}w_s(y)\,dy\right)^{p/p'}\\
 &= (C_{\alpha,J,i_0}s^\beta)^{p-1}
       \int_{D_s}|f(x(t,v,y))|^p w_s(y)\,dy .
\end{align*}
Dividing by $(C_{\alpha,J,i_0}s^\beta)^p$ gives
\[
 |F_J(t,v)|^p\le\frac{1}{C_{\alpha,J,i_0}s^\beta}
             \int_{D_s}|f(x(t,v,y))|^p w_s(y)\,dy .
\]
Since $F_J(t,v)=0$ for $s<c_dt$ and $s\simeq t$ for $s\ge c_dt$,
multiplying the above inequality by $e^{-t}t^\beta$ and integrating
in $v$ and $t$ gives, for any interval $A\subset(0,\infty)$,
\begin{align*}
 &\int_A\int |F_J(t,v)|^p\,d\omega_J(v)\,e^{-t}t^\beta\,dt
 \lesssim\int_A\int_{\{v:\,s\ge c_dt\}}\int_{D_s}
 |f(x(t,v,y))|^p e^{-t}w_s(y)\,dy\,d\omega_J(v)\,dt.
\end{align*}
Recall that $t=\Sigma(x)$, $x_j=v_j$ for $j\in J$, $x_k=y_k$ for
$k\in K$, and $x_{i_0}=s-\sum_{k\in K}y_k$. Since the change of
variables $x\mapsto(t,v,y)$ has Jacobian of absolute value one, we have
$dx=dt\,dv\,dy$. Hence
\begin{align*}
 e^{-t}w_s(y)\,dy\,d\omega_J(v)\,dt
 &=e^{-t}\left(s-\sum_{k\in K}y_k\right)^{\alpha_{i_0}}
   \prod_{k\in K}y_k^{\alpha_k}\prod_{j\in J}v_j^{\alpha_j}
   \,dy\,dv\,dt\\
 &=e^{-\Sigma(x)}\prod_{\ell=1}^d x_\ell^{\alpha_\ell}\,dx
 =d\mu_\alpha(x).
\end{align*}
Moreover, the region of integration on the right-hand side of the
preceding inequality corresponds to the set of $x\in\R^d_+$ for which
$\Sigma(x)\in A$ and
$$s=\Sigma(x)-\sum_{j\in J}x_j\ge c_d\Sigma(x)=c_dt.
$$ 
Thus, by the change of
variables, the integral on the right-hand side is equal to
\[
 \int_{\{\Sigma(x)\in A,\,s\ge c_dt\}} |f(x)|^p\,d\mu_\alpha(x)
 \le \int_{\{\Sigma(x)\in A\}} |f(x)|^p\,d\mu_\alpha(x).
\]
Consequently,
\begin{equation}\label{eq:cartesian-jensen}
 \int_A\int |F_J(t,v)|^p\,d\omega_J(v)\,e^{-t}t^\beta\,dt
 \lesssim \int_{\{\Sigma(x)\in A\}}|f(x)|^p\,d\mu_\alpha(x).
\end{equation}
If $J=\emptyset$, there are no $v$-coordinates. We write $F_\emptyset(t)$,
and the left-hand side of \eqref{eq:cartesian-jensen} becomes
\[
 \int_A |F_\emptyset(t)|^p e^{-t}t^\beta\,dt .
\]
Similarly, for every measurable function
$G\colon(0,\infty)^J\to\mathbb R$ with
$G\in L^p(d\omega_J)$, by the above change of variables
and \eqref{eq:cartesian-fiber}, we have
\begin{align*}
	 \|\one_{\{\Sigma(x)\in A,\ s\ge c_dt\}}
	G(v)\|_{L^p(d\mu_\alpha)}^p
	&=\int_A\int_{\{v:\,s\ge c_dt\}}
	|G(v)|^p e^{-t}
	\left(\int_{D_s}w_s(y)\,dy\right)d\omega_J(v)\,dt\\
	&=C_{\alpha,J,i_0}
	\int_A\int_{\{v:\,s\ge c_dt\}}
	|G(v)|^p e^{-t}s^\beta\,d\omega_J(v)\,dt\\
	&\lesssim
	\int_A e^{-t}t^\beta\,dt
	\int_{(0,\infty)^J}|G(v)|^p\,d\omega_J(v),
\end{align*}
where in the last inequality we used the fact that
$s\simeq t$ on the region $s\ge c_dt$.
It follows that
\begin{equation}\label{eq:cartesian-output}
	\|\one_{\{\Sigma(x)\in A,\ s\ge c_dt\}}
	G(v)\|_{L^p(d\mu_\alpha)}
	\lesssim
	\left(\int_A e^{-t}t^\beta\,dt\right)^{1/p}
	\|G\|_{L^p(d\omega_J)}.
\end{equation}

\subsection{Proof of Theorem~\ref{thm:main}}\label{sec:slicing}
Let $d\ge2$ and $1\le q<\infty$. By Remark~\ref{rem:reductions}, we need only to work on the admissible truncated balls. 

Let $S_{\nu}=\{x\in\R^d_+:{\nu}<\Sigma(x)\le {\nu}+1\}$ for integers ${\nu}\ge2$.
Write $f_{\nu}=f\one_{S_{\nu}}$. Recall from \eqref{eq:admissible} that the
admissible sets are
\[
 E=B_q({m},r)\cap\R^d_+,\qquad {m}\in\overline{\R^d_+},\qquad
 r\ge R_{d,q},\qquad b(E)=\inf_{z\in E}\Sigma(z)>2.
\]
Let $M$ denote the maximal operator over the admissible sets.
We prove the strip estimate with $c_p=\frac12\min\{1/p,1/p'\}$, where $p'=p/(p-1)$:
\begin{equation}\label{eq:cartesian-strip-bound}
 \|\one_{S_{\eta}}Mf_{\nu}\|_{L^p(d\mu_\alpha)}
 \lesssim e^{-c_p|{\nu}-{\eta}|}\|f_{\nu}\|_{L^p(d\mu_\alpha)},
 \qquad {\nu},{\eta}\ge2 .
\end{equation}

Once \eqref{eq:cartesian-strip-bound} has been proved, the $L^p(d\mu_\alpha)$-boundedness follows immediately {from Young's inequality for convolution on $\ell^p$.}

We now prove \eqref{eq:cartesian-strip-bound}. To do this, fix $x\in S_{\eta}$ and an admissible set $E$ containing $x$ and meeting $S_{\nu}$.
Set $b=b(E)$. If $1<q<\infty$, let $a$ be the unique minimizer of
$\Sigma$ on $\overline E$, choose $i_0$ as in \eqref{eq:i0}, and let
$L_j$ be as in \eqref{eq:Lj-strict} for $j\ne i_0$.
If $q=1$, choose $a$ as in \eqref{eq:diamond-choice}, that is,
$a=b{m}/\Sigma({m})$, choose $i_0$ as in Lemma~\ref{lem:diamond-geometry},
and let $L_j$ be as in \eqref{eq:diamond-Lj} for $j\ne i_0$.
In both cases, $\Sigma(a)=b$ and $a_{i_0}\ge b/d$.
Set
\[
 U={\nu}+1-b,\qquad V={\eta}+1-b,\qquad H=1+\max\{U,V\}.
\]
For \(y\in E\cap S_{\nu}\) and \(x\in E\cap S_{\eta}\),
\[
 2<b\le\Sigma(y)\le {\nu}+1,\qquad 2<b\le\Sigma(x)\le {\eta}+1
 \quad\Longrightarrow\quad U,V\ge0.
\] By \eqref{eq:strict-outer} for $1<q<\infty$ and
\eqref{eq:diamond-outer} for $q=1$, for a
sufficiently large fixed $C_0\ge1$, all $z\in E\cap(S_{\nu}\cup S_{\eta})$ satisfy
$|z_j-a_j|<C_0(1+H)L_j$ for $j\ne i_0$. Put
\[
 J=\{j\ne i_0:C_0(1+H)L_j<b/(2d(d-1))\}.
\]
Let \(K=\{1,\ldots,d\}\setminus(J\cup\{i_0\})\) as in
Section~\ref{sec:cartesian-marginal}, and write
\[
 I=\prod_{j\in J}\Big((a_j-C_0(1+H)L_j,\ a_j+C_0(1+H)L_j)\cap(0,\infty)\Big),
 \qquad v_x=(x_j)_{j\in J}\in I .
\]
With \({c_d}=1/(2d)\) {as in Section~\ref{sec:cartesian-marginal}}, every \(v\in I\) satisfies
\[
 \sum_{j\in J}v_j\le b-a_{i_0}+b/(2d)\le (1-{c_d})b.
\]
Recall that $s=t-\sum_{j\in J}v_j$ and $t=\Sigma(x)$ in the
coordinates of Section~\ref{sec:cartesian-marginal}. In particular,
\[
 t\ge b,\quad v\in I
 \quad\Longrightarrow\quad
 s=t-\sum_{j\in J}v_j\ge t-(1-{c_d})b\ge {c_d}t.
\]
Apply Lemma~\ref{lem:power} to each interval of \(I\). Then
\[
 \omega_J(I)\lesssim
 (1+H)^{\sum_{j\in J}(1+(\alpha_j)_+)}
 \prod_{j\in J}L_j(a_j+L_j)^{\alpha_j}.
\]
For \(k\in K\),
\[
 \frac{b}{1+H}\lesssim L_k\le a_k+L_k\lesssim b,
 \qquad
 L_k(a_k+L_k)^{\alpha_k}
 \gtrsim
 L_k^{1+(\alpha_k)_+}b^{-(-\alpha_k)_+}
 \gtrsim \frac{b^{1+\alpha_k}}{(1+H)^{1+(\alpha_k)_+}}.
\]
Recall \(\beta=\alpha_{i_0}+\sum_{k\in K}(1+\alpha_k)\), and set
\[
 D=\sum_{\ell\in J\cup K}(1+(\alpha_\ell)_+)
   =d-1+\sum_{\ell\ne i_0}(\alpha_\ell)_+ .
\]
Using Proposition~\ref{prop:mass-lower}, we obtain
\begin{equation}\label{eq:cartesian-box-size}
 \frac{\omega_J(I)}{\mu_\alpha(E)}
 \lesssim
 \frac{e^b b^{-\alpha_{i_0}}(1+H)^{\sum_{j\in J}(1+(\alpha_j)_+)}}
      {\prod_{k\in K}L_k(a_k+L_k)^{\alpha_k}}
 \lesssim e^b b^{-\beta}(1+H)^D.
\end{equation}

For \(J\ne\emptyset\), let \(\mathcal R_J\) denote the family of
all bounded open rectangles in \((0,\infty)^J\) and define
\[
 M_Jg(v)=\sup_{\substack{R\in\mathcal R_J\\v\in R}}
 \frac{1}{\omega_J(R)}\int_R |g(w)|\,d\omega_J(w).
\]
For $J=\emptyset$, set $M_\emptyset g=g$.

For each $j\in J$, let $M_j$ be the non-centered Hardy--Littlewood
maximal operator with respect to $v_j^{\alpha_j}\,dv_j$ on $(0,\infty)$,
acting on the $j$-th coordinate. By \eqref{eq:power-centered}, this
measure is doubling. Hence,
\[
 \|M_jg\|_{L^p(v_j^{\alpha_j}\,dv_j)\to L^p(v_j^{\alpha_j}\,dv_j)}\lesssim 1,
 \qquad 1<p<\infty.
\]
Writing $J=\{j_1,\ldots,j_k\}$, we have
\[
 M_Jg\le M_{j_1}\cdots M_{j_k}g.
\]
Consequently,
\[
 \|M_Jg\|_{L^p(d\omega_J)}\lesssim \|g\|_{L^p(d\omega_J)},
 \qquad 1<p<\infty.
\]
The last estimate also holds for $J=\emptyset$.

For each pair $(i_0,J)$ with $i_0\notin J$, write $F_{i_0,J}$ for
the function $F_J$ defined in Section~\ref{sec:cartesian-marginal}
with this choice of $i_0$ and $J$.
For \(\max\{{\nu},b\}<t\le {\nu}+1\) and \(v\in I\), take
\(s=t-\sum_{j\in J}v_j\) and put
\[
 E_{t,v}=\{y\in D_s:x(t,v,y)\in E\}.
\]
For fixed \(v,y\), the substitution
\(x_{i_0}=t-\sum_{j\in J}v_j-\sum_{k\in K}y_k\) gives
\(dx_{i_0}=dt\). Since \((x_j)_{j\in J}\in I\) on \(E\cap S_{\nu}\),
\[
 \int_E |f_{\nu}|\,d\mu_\alpha
 =\int_{\max\{{\nu},b\}}^{{\nu}+1}\int_I e^{-t}
   \int_{E_{t,v}} |f(x(t,v,y))|w_s(y)\,dy\,d\omega_J(v)\,dt.
\]
Since \(E_{t,v}\subseteq D_s\) and \(s\ge {c_d}t\), the definition of \(F_{i_0,J}\)
gives
\[
 \int_{E_{t,v}} |f(x(t,v,y))|w_s(y)\,dy
 \le \int_{D_s} |f(x(t,v,y))|w_s(y)\,dy
 = C_{\alpha,J,i_0}s^\beta F_{i_0,J}(t,v).
\]
Substituting this bound and using $v_x\in I$ gives
\begin{align*}
 \frac{1}{\mu_\alpha(E)}\int_E |f_{\nu}|\,d\mu_\alpha
 &\le \frac{C_{\alpha,J,i_0}}{\mu_\alpha(E)}
 \int_{\max\{{\nu},b\}}^{{\nu}+1}\int_I
 F_{i_0,J}(t,v)e^{-t}s^\beta\,d\omega_J(v)\,dt\\
 &\lesssim \frac{\omega_J(I)}{\mu_\alpha(E)}
 \int_{\nu}^{{\nu}+1} M_J(F_{i_0,J}(t,\cdot))(v_x)e^{-t}t^\beta\,dt\\
 &\lesssim e^b b^{-\beta}(1+H)^D
 \int_{\nu}^{{\nu}+1} M_J(F_{i_0,J}(t,\cdot))(v_x)e^{-t}t^\beta\,dt,
\end{align*}
where in the second inequality we used ${c_d}t\le s\le t$.

For fixed \((i_0,J)\) and input strip \(S_{\nu}\), define
\begin{align*}
 \sigma_\ell&=\int_\ell^{\ell+1}e^{-t}t^\beta\,dt
 \simeq e^{-\ell}(\ell+1)^\beta, \quad \ell\ge2,\\
 G_{i_0,J}(v)&=\left(\int_{\nu}^{{\nu}+1}[M_J(F_{i_0,J}(t,\cdot))(v)]^p
                 e^{-t}t^\beta\,dt\right)^{1/p}.
\end{align*}
Using ${\nu}+1=b+U$ and ${\eta}+1=b+V$, we obtain
\begin{align*}
 e^b b^{-\beta}(1+H)^D \sigma_{\nu}^{1/p'}\sigma_{\eta}^{1/p}
 &\lesssim (1+H)^D e^{-U/p'-V/p}
      (1+U/b)^{\beta/p'}(1+V/b)^{\beta/p}\\
 &\lesssim (1+H)^{D+\beta_+}
       e^{-\min\{1/p',1/p\}(H-1)}
 \lesssim_{d,p,q,\alpha}e^{-c_p|{\nu}-{\eta}|},
\end{align*}
since $H-1=\max\{U,V\}\ge |{\nu}-{\eta}|$.
By H\"older's inequality,
\begin{align*}
 &\int_{\nu}^{{\nu}+1}M_J(F_{i_0,J}(t,\cdot))(v_x)e^{-t}t^\beta\,dt\\
 &\quad\le
 \left(\int_{\nu}^{{\nu}+1}e^{-t}t^\beta\,dt\right)^{1/p'}
 \left(\int_{\nu}^{{\nu}+1}[M_J(F_{i_0,J}(t,\cdot))(v_x)]^p
 e^{-t}t^\beta\,dt\right)^{1/p}\\
 &\quad=\sigma_{\nu}^{1/p'}G_{i_0,J}(v_x).
\end{align*}
Combining the preceding estimates, we obtain
\begin{align*}
 \frac{1}{\mu_\alpha(E)}\int_E |f_{\nu}|\,d\mu_\alpha
 &\lesssim e^b b^{-\beta}(1+H)^D\sigma_{\nu}^{1/p'}G_{i_0,J}(v_x)\\
 &=\bigl[e^b b^{-\beta}(1+H)^D\sigma_{\nu}^{1/p'}\sigma_{\eta}^{1/p}\bigr]
 \sigma_{\eta}^{-1/p}G_{i_0,J}(v_x)\\
 &\lesssim e^{-c_p|{\nu}-{\eta}|}\sigma_{\eta}^{-1/p}G_{i_0,J}(v_x).
\end{align*}
Since $v_x\in I$ and $\Sigma(x)\ge b$, we also have
\[
 \Sigma(x)-\sum_{j\in J}x_j
 \ge\Sigma(x)-(1-{c_d})b\ge {c_d}\Sigma(x).
\]
Set
\[
 A_{{\eta},J}=\left\{x\in S_{\eta}:\Sigma(x)-\sum_{j\in J}x_j
 \ge {c_d}\Sigma(x)\right\}.
\]
Thus $x\in A_{{\eta},J}$, and
\[
 \frac{1}{\mu_\alpha(E)}\int_E|f_{\nu}|\,d\mu_\alpha
 \lesssim e^{-c_p|{\nu}-{\eta}|}\sigma_{\eta}^{-1/p}
 \one_{A_{{\eta},J}}(x)G_{i_0,J}(v_x).
\]
The implicit constants are independent of $E$, ${\nu}$, ${\eta}$, $i_0$ and $J$.
Since $E$ is arbitrary and admissible sets not meeting $S_{\nu}$ contribute
zero, we obtain
\begin{align*}
 \one_{S_{\eta}}(x)Mf_{\nu}(x)
 &\lesssim e^{-c_p|{\nu}-{\eta}|}
 \max_{i_0,J}\left\{\sigma_{\eta}^{-1/p}\one_{A_{{\eta},J}}(x)
 G_{i_0,J}((x_j)_{j\in J})\right\}\\
 &\le e^{-c_p|{\nu}-{\eta}|}
 \sum_{i_0,J}\sigma_{\eta}^{-1/p}\one_{A_{{\eta},J}}(x)
 G_{i_0,J}((x_j)_{j\in J}).
\end{align*}
Here the maximum and sum are over all pairs $i_0\in\{1,\ldots,d\}$ and
$J\subseteq\{1,\ldots,d\}\setminus\{i_0\}$, and $\sigma_{\eta}$ in each
term is computed with the corresponding $\beta$.
Thus \eqref{eq:cartesian-output} gives
\[
 \|\one_{S_{\eta}}Mf_{\nu}\|_{L^p(d\mu_\alpha)}
 \lesssim e^{-c_p|{\nu}-{\eta}|}\sum_{i_0,J}\|G_{i_0,J}\|_{L^p(d\omega_J)}.
\]
Finally, Fubini's theorem, the $L^p(d\omega_J)$ bound for $M_J$ and
\eqref{eq:cartesian-jensen} give
\begin{align*}
 \|G_{i_0,J}\|_{L^p(d\omega_J)}^p
 &=\int_{\nu}^{{\nu}+1}\|M_J(F_{i_0,J}(t,\cdot))\|_{L^p(d\omega_J)}^p
               e^{-t}t^\beta\,dt\\
 &\lesssim\int_{\nu}^{{\nu}+1}\|F_{i_0,J}(t,\cdot)\|_{L^p(d\omega_J)}^p
               e^{-t}t^\beta\,dt
 \lesssim\|f_{\nu}\|_{L^p(d\mu_\alpha)}^p .
\end{align*}
Since there are only $d2^{d-1}$ possible pairs $(i_0,J)$, this proves
\eqref{eq:cartesian-strip-bound}.

This completes our proof.

\section{\texorpdfstring{Failure of weak type $(1,1)$}{Failure of weak type (1,1)}}\label{sec:weak-endpoint}

\subsection{The strictly convex range \texorpdfstring{$1<q<\infty$}{1 < q < infinity}}
\begin{proposition}\label{prop:weak-strict}
For every $d\ge2$, $1<q<\infty$ and $\alpha\in(-1,\infty)^d$,
$M_{\mu_\alpha,q}$ is not of weak type $(1,1)$.
More precisely, there exist $f_s\ge0$ with $\|f_s\|_{L^1(d\mu_\alpha)}=1$
and $\lambda_s>0$, defined for all sufficiently large \(s\), such that
\begin{equation}\label{eq:weak-strict-growth}
 \lambda_s\mu_\alpha\bigl(\{M_{\mu_\alpha,q}f_s>\lambda_s\}\bigr)
 \gtrsim_{d,q,\alpha}s^{(d-1)/2}\longrightarrow\infty
 \qquad (s\to\infty).
\end{equation}
\end{proposition}
\begin{proof}
Set $\mathbf e=(1,\ldots,1)$, $u_q=d^{-1/q}\mathbf e$, and
\[
 H=\{y\in\R^d:\Sigma(y)=0\},\qquad
 Y_s=\{y\in H:|y|_q<s/8\},\qquad x(s)=2s\mathbf e.
\]

Let $s$ be sufficiently large. For $y\in Y_s$, put

\[
 B_y=B_q(x(s)+y,s),\qquad D_s=B_q(x(s),1).
\]

For $x\in B_y$, $7s/8<x_j<25s/8$. By Lemma~\ref{lem:lowest}, the minimum of $\Sigma$ on $\overline{B_y}$ is attained at

\[
 a_y=x(s)+y-su_q,
 \qquad b_s=\Sigma(a_y)=2ds-sd^{1/q'},
\]
where $q'$ is conjugate to $q$. 
Note that $(a_y)_j\simeq s$ for every $j$ and hence $L_j=\sqrt s$ in \eqref{eq:Lj-strict}.

By Lemma~\ref{lem:strict-geometry}(ii) and $L_j=\sqrt s$,
\[
 \mu_0(B_y)\gtrsim_{d,q}\int_1^2e^{-b_s-h}s^{(d-1)/2}\,dh\simeq e^{-b_s}s^{(d-1)/2}.
\]
By Lemma~\ref{lem:strict-geometry}(i), the {$(d-1)$-dimensional Lebesgue measure of each section in the coordinates $(x_j)_{j\ne i_0}$} is at most $C_{d,q}(1+h)^{d-1}s^{(d-1)/2}$. Thus
\[
 \mu_0(B_y)\lesssim_{d,q}e^{-b_s}s^{(d-1)/2}\int_0^\infty e^{-h}(1+h)^{d-1}\,dh
 \lesssim_{d,q}e^{-b_s}s^{(d-1)/2}.
\]
{Since $x_j\simeq s$ for $x\in B_y$ and $1\le j\le d$, we obtain}
\begin{equation}\label{eq:weak-strict-mass}
 \mu_\alpha(B_y)\simeq_{d,q,\alpha}s^{\sum_j\alpha_j+(d-1)/2}e^{-b_s},\qquad y\in Y_s.
\end{equation}

For sufficiently large $s$, $D_s\subset B_y$ for every $y\in Y_s$,
because $1+s/8<s$. Take
$f_s=\mu_\alpha(D_s)^{-1}\one_{D_s}$. Define
\[
 \sigma_s=2s-sd^{-1/q},\qquad
 \mathcal E_s=\{\sigma\mathbf e+y:\sigma_s+2<\sigma<\sigma_s+3,
                       \,y\in Y_s\}.
\]
If $x=\sigma\mathbf e+y\in\mathcal E_s$, then for large $s$ we have
$\sigma<2s$ and
\[
 |x-(x(s)+y)|_q=d^{1/q}(2s-\sigma)<s.
\]
Thus $x\in B_y$. The inclusion $D_s\subset B_y$ and
\eqref{eq:weak-strict-mass} imply
\begin{equation}\label{eq:weak-strict-pointwise}
 M_{\mu_\alpha,q}f_s(x)
 \ge \mu_\alpha(B_y)^{-1}
 \gtrsim_{d,q,\alpha} e^{b_s}s^{-\sum_j\alpha_j-(d-1)/2},
 \qquad x\in\mathcal E_s.
\end{equation}

Since $\Sigma(\sigma\mathbf e+y)=d\sigma$, $d\sigma_s=b_s$ and $x_j\simeq s$ on $\mathcal E_s$, we have
\begin{equation}\label{eq:weak-strict-output}
 \mu_\alpha(\mathcal E_s)\simeq s^{\sum_j\alpha_j}\int_{\sigma_s+2}^{\sigma_s+3}e^{-d\sigma}\,d\sigma\int_{Y_s}dy
 \simeq_{d,q,\alpha}e^{-b_s}s^{\sum_j\alpha_j+d-1}.
\end{equation}
{Here $dy$ denotes Lebesgue measure on the hyperplane $H$.}
Set $\lambda_s=c e^{b_s}s^{-\sum_j\alpha_j-(d-1)/2}$ with $c>0$ sufficiently small. By \eqref{eq:weak-strict-pointwise},
\[
 \lambda_s\mu_\alpha\{M_{\mu_\alpha,q}f_s>\lambda_s\}
 \ge\lambda_s\mu_\alpha(\mathcal E_s)\gtrsim s^{(d-1)/2}.
\]

{This completes our proof.}
\end{proof}

\subsection{Cubes, \texorpdfstring{$q=\infty$}{q = infinity}}

Failure of weak type \((1,1)\) for cubes with respect to \(\mu_\alpha\)
was already proved in \cite[Remark~3.1]{NSSS}. Here, we refine their argument to give a simpler proof.
\begin{proposition}\label{prop:weak-cubes}
If $d\ge2$ and $\alpha\in(-1,\infty)^d$, then
$M_{\mu_\alpha,\infty}$ is not of weak type $(1,1)$.
\end{proposition}
\begin{proof}

Let $x(s)=2s\mathbf e$ and $H$ be as in the proof of Proposition~\ref{prop:weak-strict}. Put
\[
 Y_s=\{y\in H:|y|_\infty<s/8\},\quad Q_y=B_\infty(x(s)+y,s),\quad D_s=B_\infty(x(s),1).
\]
For $s$ large, $D_s\subset Q_y\subset\R^d_+$ for all $y\in Y_s$. Moreover,
\begin{align*}
 \mu_\alpha(Q_y)&=\prod_{j=1}^d\int_{s+y_j}^{3s+y_j}t^{\alpha_j}e^{-t}\,dt
 \simeq_{d,\alpha}s^{\sum_j\alpha_j}e^{-ds-\Sigma(y)}=s^{\sum_j\alpha_j}e^{-ds}.
\end{align*}
Let $\mathcal E_s=\{\sigma\mathbf e+y:s+2<\sigma<s+3,\ y\in Y_s\}$. For $x=\sigma\mathbf e+y\in\mathcal E_s$, {when $s$ is sufficiently large,} $|x-(x(s)+y)|_\infty=2s-\sigma<s$. Hence, for $f_s=\mu_\alpha(D_s)^{-1}\one_{D_s}$,
\[
 M_{\mu_\alpha,\infty}f_s(x)\ge\frac1{\mu_\alpha(Q_y)}\gtrsim e^{ds}s^{-\sum_j\alpha_j},\qquad x\in\mathcal E_s.
\]
As in the proof of Proposition~\ref{prop:weak-strict},
\[
 \mu_\alpha(\mathcal E_s)\simeq_{d,\alpha}s^{\sum_j\alpha_j+d-1}e^{-ds}.
\]
For $\lambda_s=ce^{ds}s^{-\sum_j\alpha_j}$ with $c>0$ sufficiently small,
\[
 \lambda_s\mu_\alpha\{M_{\mu_\alpha,\infty}f_s>\lambda_s\}\ge\lambda_s\mu_\alpha(\mathcal E_s)\gtrsim s^{d-1}\to\infty.
\]

{This completes our proof.}
\end{proof}

\subsection{Diamonds, \texorpdfstring{$q=1$}{q = 1}}

We shall use the following estimate.

\begin{lemma}
\label{lem:diamond-upper-weak}
Let $d\ge2$ and $\alpha_j>-1$. Take $\xi=(\xi',\xi_d)\in\R^d_+$,
where \(\xi'=(\xi_1,\ldots,\xi_{d-1})\), satisfying
\begin{equation}\label{eq:weak-diamond-conditions}
 \xi_d\ge2,\qquad\sum_{j<d}\xi_j\le\xi_d.
\end{equation}
Let
\[
 D_{\xi,N}=B_1\bigl((\xi',\xi_d+3N),3N+1\bigr),\qquad N>0.
\]
Then
\begin{equation}\label{eq:weak-diamond-mass}
 \mu_\alpha(D_{\xi,N})\le C_{d,\alpha}
 e^{-\Sigma(\xi)}\xi_d^{\alpha_d}
 \prod_{j<d}(1+\xi_j)^{\alpha_j+1},
\end{equation}
with a constant independent of $\xi$ and $N$.
\end{lemma}
\begin{proof}
For $x\in D_{\xi,N}\cap\R^d_+$, the diamond inequality implies
\[
 x_d>\xi_d+\sum_{j<d}|x_j-\xi_j|-1.
\]
Set
\[
 K_-=\sum_{j<d}(\xi_j-x_j)_+,
 \qquad K_+=\sum_{j<d}(x_j-\xi_j)_+.
\]
For fixed \(x'=(x_1,\ldots,x_{d-1})\in(0,\infty)^{d-1}\), set
\[
 u=\xi_d+K_-+K_+-1,\qquad 0\le K_-\le\xi_d,\qquad
 \xi_d/2\le u\le 2\xi_d+K_+ .
\]
For \(u\ge1\),
\begin{equation}\label{eq:weak-gamma-tail}
 \int_u^\infty t^{\alpha_d}e^{-t}\,dt
 \le C_{\alpha_d}u^{\alpha_d}e^{-u}.
\end{equation}
The lower bound for \(u\) gives the first case below, and the upper
bound gives the second:
\[
 u^{\alpha_d}\le
 \begin{cases}
  2^{-\alpha_d}\xi_d^{\alpha_d},&\alpha_d<0,\\
  C_{\alpha_d}\xi_d^{\alpha_d}(1+K_+)^{\alpha_d}
  \le C_{\alpha_d}\xi_d^{\alpha_d}e^{K_+/2},&\alpha_d\ge0 .
 \end{cases}
\]
Since \(e^{-u}=e\,e^{-\xi_d}e^{-K_-}e^{-K_+}\), we obtain
\begin{equation}\label{eq:weak-diamond-tail}
 \int_u^\infty t^{\alpha_d}e^{-t}\,dt
 \le C_{\alpha_d}\xi_d^{\alpha_d}e^{-\xi_d}e^{-K_-}e^{-K_+/2}.
\end{equation}
For this \(x'\), the diamond inequality gives
\[
 \{x_d>0:(x',x_d)\in D_{\xi,N}\}\subset(u,\infty).
\]
Hence
\[
 \mu_\alpha(D_{\xi,N}\cap\R^d_+)
 \le \int_{(0,\infty)^{d-1}}
 \left(\prod_{j<d} x_j^{\alpha_j}e^{-x_j}\right)
 \left(\int_u^\infty t^{\alpha_d}e^{-t}\,dt\right)\,dx'.
\]
Inserting \eqref{eq:weak-diamond-tail}, the factors separate since
\[
 e^{-K_-}e^{-K_+/2}
 =\prod_{j<d}e^{-(\xi_j-x_j)_+}e^{-(x_j-\xi_j)_+/2}.
\]
Now \(x_j+(\xi_j-x_j)_+=\xi_j+(x_j-\xi_j)_+\), so
\[
 \mu_\alpha(D_{\xi,N}\cap\R^d_+)
 \lesssim_{\alpha_d} e^{-\Sigma(\xi)}\xi_d^{\alpha_d}
 \prod_{j<d}\int_0^\infty x_j^{\alpha_j}
 e^{-\frac32(x_j-\xi_j)_+}\,dx_j.
\]
For $\beta>-1$ and $a\ge0$, we have
\begin{align}
 \int_0^\infty t^\beta e^{-\frac32(t-a)_+}\,dt
 &=\frac{a^{\beta+1}}{\beta+1}+\int_0^\infty(a+z)^\beta e^{-3z/2}\,dz\notag\\
 &\lesssim_\beta (1+a)^{\beta+1}.\label{eq:weak-diamond-transverse}
\end{align}
Here we used
\[
 \int_0^\infty(a+z)^\beta e^{-3z/2}\,dz\lesssim_\beta
 \begin{cases}
 1,&0\le a\le1,\\
 a^\beta,&a>1.
 \end{cases}
\]
By \eqref{eq:weak-diamond-transverse},
\[
 \mu_\alpha(D_{\xi,N}\cap\R^d_+)\lesssim_{d,\alpha}e^{-\Sigma(\xi)}\xi_d^{\alpha_d}\prod_{j<d}(1+\xi_j)^{\alpha_j+1}.
\]

{
By \eqref{eq:fold}, since the center of $D_{\xi,N}$ lies in $\overline{\R^d_+}$, we have
\[
 T(D_{\xi,N})\cap\R^d_+=D_{\xi,N}\cap\R^d_+.
\]
By \eqref{eq:weight},
\begin{align*}
 \mu_\alpha(D_{\xi,N})
 &\le 2^d\mu_\alpha(T(D_{\xi,N})\cap\R^d_+)\\
 &=2^d\mu_\alpha(D_{\xi,N}\cap\R^d_+).
\end{align*}
}
This completes our proof.
\end{proof}

\begin{proposition}\label{prop:weak-diamonds}
Let $d\ge2$ and $\alpha\in(-1,\infty)^d$. Then
$M_{\mu_\alpha,1}$ is not of weak type $(1,1)$.
More precisely, there are $f_N\ge0$ with
$\|f_N\|_{L^1(d\mu_\alpha)}=1$ and $\lambda_N>0$, defined for all
sufficiently large \(N\), such that
\begin{equation}\label{eq:weak-diamond-log}
 \lambda_N\mu_\alpha\bigl(\{M_{\mu_\alpha,1}f_N>
 \lambda_N\}\bigr)\gtrsim_{d,\alpha}\log N
 \longrightarrow\infty\qquad (N\to\infty).
\end{equation}
\end{proposition}
\begin{proof}
Put
\begin{equation}\label{eq:theta}
 \Theta=\sum_{j<d}(\alpha_j+1)>0.
\end{equation}
Fix $0<\eta<1/(4d)$. For $N>0$, set
\[
 F_N=(0,\eta)^{d-1}\times(N,N+\eta),\qquad
 f_N=\frac{\one_{F_N}}{\mu_\alpha(F_N)}
\]
so that $\|f_N\|_{L^1(d\mu_\alpha)}=1$.

Fix \(\varepsilon=1/[16(d-1)]\), and take \(N\) so large that
\[
 2<\tfrac{\Theta}{2}\log N\le\tfrac{N}{8},
 \qquad \varepsilon\sqrt N\ge1.
\]
For \(2<v<\frac{\Theta}{2}\log N\), set \(L_N(v)=\varepsilon N e^{-v/\Theta}\).
Let \(\mathcal E_N\) consist of the points \(\xi\) given by
\begin{equation}\label{eq:diamond-output-param}
 \begin{split}
 &0<\xi_j<L_N(v),\qquad 1\le j<d,\\
 &\xi_d=N-v-\sum_{j<d}\xi_j,
 \qquad 2<v<\tfrac{\Theta}{2}\log N.
 \end{split}
\end{equation}
{For such points,} \(v\le N/8\) and \(\sum_{j<d}\xi_j\le(d-1)\varepsilon N<N/8\).
Also \(N^{-1/2}\le e^{-v/\Theta}\le1\). Thus on \(\mathcal E_N\),
\begin{equation}\label{eq:diamond-output-conditions}
 \tfrac34N\le\xi_d\le N,\qquad
 \sum_{j<d}\xi_j\le\tfrac18N\le\xi_d,\qquad
 1\le\varepsilon\sqrt N\le L_N(v)\le\varepsilon N.
\end{equation}
In particular, the hypotheses of Lemma~\ref{lem:diamond-upper-weak}
are satisfied.

For $\xi\in\mathcal E_N$, consider the full-space diamond
\[
 D_{\xi,N}=B_1((\xi',\xi_d+3N),3N+1).
\]
It contains $\xi$, since the distance from $\xi$ to its center is
$3N<3N+1$. It also contains all of $F_N$. Indeed, for $y\in F_N$,
$y_d<\xi_d+3N$ and
\begin{align*}
 |y-(\xi',\xi_d+3N)|_1
 &=\sum_{j<d}|y_j-\xi_j|+\xi_d+3N-y_d\\
 &\le\sum_{j<d}(\xi_j+y_j)+\xi_d+3N-N\\
 &\le3N-v+(d-1)\eta<3N+1,
\end{align*}
where we used $v>2$ and $(d-1)\eta<1$.
Therefore
\begin{equation}\label{eq:diamond-M-lower}
 M_{\mu_\alpha,1}f_N(\xi)
 \ge\frac1{\mu_\alpha(D_{\xi,N})},
 \qquad\xi\in\mathcal E_N.
\end{equation}

By Lemma~\ref{lem:diamond-upper-weak} and
\eqref{eq:diamond-output-conditions},
\begin{align*}
 \mu_\alpha(D_{\xi,N})
 &\lesssim_{d,\alpha}
 e^{-N+v}N^{\alpha_d}\prod_{j<d}(1+L_N(v))^{\alpha_j+1}\\
 &\lesssim_{d,\alpha}
 e^{-N+v}N^{\alpha_d} L_N(v)^\Theta
 {=\varepsilon^\Theta e^{-N}N^{\alpha_d+\Theta}},
\end{align*}
where we used $L_N(v)\ge1$ and
$e^vL_N(v)^\Theta=\varepsilon^\Theta N^\Theta$.
It follows that
\begin{equation}\label{eq:diamond-pointwise-lower}
 M_{\mu_\alpha,1}f_N(\xi)
 \ge c_{d,\alpha}e^NN^{-\alpha_d-\Theta},
 \qquad\xi\in\mathcal E_N.
\end{equation}
Choose
\begin{equation}\label{eq:diamond-lambda}
 \lambda_N=\tfrac12c_{d,\alpha}e^NN^{-\alpha_d-\Theta},
\end{equation}
so that $\mathcal E_N\subset\{M_{\mu_\alpha,1}f_N>\lambda_N\}$.

The parametrization \eqref{eq:diamond-output-param} has Jacobian of
absolute value one. Since $\Sigma(\xi)=N-v$ and
$\xi_d\simeq N$, we obtain
\begin{align}
 \mu_\alpha(\mathcal E_N)
 &\gtrsim_{d,\alpha}
 N^{\alpha_d}e^{-N}\int_{2}^{(\Theta/2)\log N}e^v
 \prod_{j<d}\left(\int_0^{L_N(v)}t^{\alpha_j}\,dt\right)dv
 \notag\\
 &=\frac{N^{\alpha_d}e^{-N}}{\prod_{j<d}(\alpha_j+1)}
 \int_2^{(\Theta/2)\log N}e^vL_N(v)^\Theta\,dv
 \notag\\
 &=\frac{\varepsilon^\Theta N^{\alpha_d+\Theta}e^{-N}}{\prod_{j<d}(\alpha_j+1)}
 \left(\tfrac{\Theta}{2}\log N-2\right)
 \gtrsim_{d,\alpha}e^{-N}N^{\alpha_d+\Theta}\log N.
 \label{eq:diamond-output-mass}
\end{align}
Together with \eqref{eq:diamond-lambda}, this proves
\eqref{eq:weak-diamond-log}.

{This completes our proof.}
\end{proof}

\begin{proof}[Proof of  Theorem~\ref{thm:endpoint}]
For $d\ge2$, Propositions~\ref{prop:weak-strict},
\ref{prop:weak-cubes}, and \ref{prop:weak-diamonds} give failure of
weak type $(1,1)$ for $1<q<\infty$, $q=\infty$, and $q=1$,
respectively. For $d=1$, every $\ell^q$ ball is an interval;
the non-centered interval maximal operator is of weak type $(1,1)$
with respect to every locally finite Borel measure on the line.

{This completes our proof.}
\end{proof}

\noindent{\bf AI use declaration.}
During the preparation of this manuscript, the authors used ChatGPT 6.0 to assist with language editing, the organisation and drafting of parts of the exposition, and the review and reformulation of selected mathematical arguments and estimates.  All AI-assisted material was reviewed and edited by the authors, who take full responsibility for the accuracy, originality, proofs, references, and conclusions of the manuscript.

\medskip
\noindent{\bf Conflict of interest.}
The authors declare that they have no conflict of interest.

\medskip
\noindent{\bf Data availability statement.} No datasets were generated or analysed, so data sharing is not applicable.

\medskip
\noindent{\bf Acknowledgements.}
The first author was supported by the Australian Research Council Discovery Project DP260101083.

\end{document}